\documentclass[]{article}

\usepackage[intlimits]{amsmath}
\usepackage{natbib}
\usepackage{amsthm,amssymb,amsfonts}
\usepackage{hyperref}

\bibpunct{(}{)}{,}{}{}{,}

\newtheorem{thm}{Theorem}[section]      

\newcommand{\BT}{\begin{thm}}   \newcommand{\ET}{\end{thm}}

\newtheorem{lem}[thm]{Lemma}
\newcommand{\BL}{\begin{lem}}   \newcommand{\EL}{\end{lem}}
\newtheorem{prop}[thm]{Proposition}
\newcommand{\BP}{\begin{prop}}   \newcommand{\EP}{\end{prop}}
\newtheorem{clm}[thm]{Claim}            %
\newcommand{\BCM}{\begin{samepage} \begin{clm}}   \newcommand{\ECM}{\end{clm} \end{samepage}}

\theoremstyle{definition}
\newtheorem{remark}[thm]{Remark}
\newcommand{\BRM}{\begin{remark}} \newcommand{\ERM}{\end{remark}}
\newtheorem{dfn}[thm]{Definition}      %
\newcommand{\BD}{\begin{dfn}}   \newcommand{\ED}{\end{dfn}}
\newtheorem{construction}[thm]{Construction}
\newcommand{\BCON}{\begin{construction}} \newcommand{\ECON}{\end{construction}}
\newtheorem{notations}[thm]{Notation}
\newcommand{\BNOT}{\begin{notations}} \newcommand{\ENOT}{\end{notations}}
\newtheorem{exmp}[thm]{Example}
\newcommand{\BEX}{\begin{exmp}} \newcommand{\EEX}{\end{exmp}}

\newtheorem{fact}[thm]{Fact}            %
\newcommand{\BF}{\begin{fact}}   \newcommand{\EF}{\end{fact}}
\newtheorem{techcorr}[thm]{Corollary}      %
\newcommand{\Bcr}{\begin{techcorr}}
\newcommand{\Ecr}{\end{techcorr}}
\newcommand{\BE}{\begin{enumerate}}
\newcommand{\EE}{\end{enumerate}}
\newcommand{\BI}{\begin{itemize}}
\newcommand{\EI}{\end{itemize}}
\newcommand{\BPF}{\begin{proof}} \newcommand {\EPF}{\end{proof}}

\newcommand{\R}{{\mathbb{R}}}
\newcommand{\N}{{\mathbb{N}}}

\newcommand{\Z}{{\mathbb{Z}}}

\DeclareMathOperator{\adj}{adj} %
\DeclareMathOperator{\sgn}{sgn} %
\DeclareMathOperator{\spn}{span} %

\numberwithin{equation}{section}
\newcommand{\system}{$(A_k,b_k,c_k)_{k \in \Z}$ }
\newcommand{\tildesystem}{$(\tilde{A}_k,\tilde{b}_k,\tilde{c}_k)_{k \in \Z}$ }
\newcommand{\REF}[2]{\hyperref[#2]{#1~\ref*{#2}}}

\begin{document}

\title{Memoryless output nullification and canonical forms, for time varying systems}

\author{Gera Weiss\thanks{Department of Computer Science and Applied Mathematics, The Weizmann Institute of Science,
Rehovot, 76100, Israel. gera.weiss@weizmann.ac.il.} \thanks{Research supported by grants from the Israel Science Foundation and from the Information
Society Technologies Programme of the European Commission.}}

\date{March 15, 2005}

\maketitle

\begin{abstract}
We study the possibility of nullifying time-varying systems with memoryless output feedback. The systems we examine are linear single-input
single-output finite-dimensional time-varying systems. For generic completely controllable and completely observable discrete-time systems, we show that
any state at any time can be steered to the origin within finite time. An algorithm for nullification and an upper bound for nullification time,
depending only on the system's dimension, are provided. The algorithm is described using a representation of the system in time-varying controller
canonical form. We verify that every completely controllable system has such a representation. The application of the nullification algorithm to
sampled-data systems is also analysed: we show that a controllable continuous-time time-varying system with analytic coefficients can be nullified
utilising zero-hold sampling of the output and time-varying memoryless linear feedback; for generic observables, almost any sampling period can be used.
We also prove that controllability of time-varying systems with analytic coefficients is preserved under zero-hold sampling at almost any rate.
\end{abstract}

\section{Introduction and main results}

As a sharp form of stabilisation, state nullification is appealing both from a mathematical perspective and from the applications point of view. This
paper examines the possibility of nullifying time-varying systems with memoryless output feedback. While efficient algorithms for stabilisation invoking
more complex dynamic feedback are available, the use of memoryless feedback strategy has apparent mathematical and engineering advantages due to its
simplicity. The time invariant case was analysed in~\citet{AW04}. The extension to time-varying systems, provided here, employs two new tools, namely,
controller canonical form and preservation of controllability under sampling.


\medskip

Four main results are presented: two concerning discrete-time systems and two concerning continuous-time systems.

\medskip

The main results about discrete-time systems are as follows.

Consider a finite dimensional linear time-varying discrete-time scalar-input scalar-output control systems of the form
\begin{equation} \label{equ:system}
\begin{aligned} x_{k+1} &= A_k x_k + b_k u_k \\
y_k &= c_k x_k
\end{aligned}\end{equation}
where, for every $k \in \Z$,  $A_k$ is $n \times n$ matrix, $b_k$ is $n$~dimensional column vector and $c_k$ is $n$~dimensional row vector. The data
specifying a concrete system is given by the doubly infinite sequence \system\!\!.

In most of this text, attention is restricted to completely controllable and completely observable systems as given by
the following definitions.

\BD The control system~\eqref{equ:system} is \emph{completely controllable} if for every $k \in \Z$ and every pair of states, $\xi_s, \xi_f \in \R^n$,
there are controls $u_k,...,u_{k+n-1}$ such that if $x_k = \xi_s$ then $x_{k+n} = \xi_f$. \ED

\BD The control system~\eqref{equ:system} is \emph{completely observable} if for every $k \in \Z$ and controls $u_k,...,u_{k+n-2}$, the state $x_k$ is
determined uniquely by the observations $y_k,...,y_{k+n-1}$. \ED

The control strategy that we consider is time-varying memoryless feedback from the output in the form
\begin{equation} \label{equ:controller}
u_k = F_k y_k
\end{equation}
where $\{F_k\}_{k \in \Z}$ is a sequence of scalars. When such a feedback is applied, the dynamics have the form
\begin{equation} \label{equ:dynamics}
x_{k+1} = (A_k + F_k b_k c_k) x_k.
\end{equation}

\medskip

The control objective studied in this paper is state nullification of time-varying discrete-time systems by memoryless output feedback, as given in the
following definition.

\BD \label{def:nullification} The system~\eqref{equ:system} is \emph{uniformly nullifiable by memoryless linear output feedback} if there is a sequence
of scalars $\{F_k\}_{k \in \Z}$ and a constant $N \in \N$ such that, for every $k \in \Z$ and every $x_k \in \R^n$, the sequence $x_k,
x_{k+1},...,x_{k+N}$ resulting from the dynamics~\eqref{equ:dynamics} satisfies $x_{k+N}=0$. \ED

This objective is related to~\citet{Bro99}, where the following open problem is offered: find a linear memoryless output feedback such that the
resulting closed-loop system is uniformly exponentially stable. The original problem is stated for continuous-time time-invariant systems. Discrete-time
analogues are studied in e.g.~\citep{AeyelsWillems92,Leonov02,AW04}.

For the formulation of the following theorem, recall the notion of the adjugate, or adjoint, of an $n \times n$ matrix. The adjugate matrix is denoted
by $\adj(A)$. It is the $n \times n$ matrix whose entry in row $j$ and column $i$ is given by $(-1)^{i+j}M_{ij}$ where $M_{ij}$ represents the $(n-1)
\times (n-1)$ minor of $A$ obtained by deleting row $i$ and column $j$ of $A$ and taking the determinant of the resulting $(n-1) \times(n-1)$ matrix
\citep[pages 56, 85]{Hoh64}.

\medskip

The first result concerning discrete-time systems is as follows.

\BT \label{thm:nullification} If the system~\eqref{equ:system} is completely controllable, completely observable and, for every $k \in \Z$, $c_k
\adj(A_k) b_k \ne 0$, then it is uniformly nullifiable by memoryless linear output feedback.  \ET

Our proof of \REF{Theorem}{thm:nullification} gives an algorithm that finds feedback coefficients for nullification. For that algorithm, we show that
the number of steps needed for nullification (the number $N$ mentioned in \REF{Definition}{def:nullification}) is bounded by $2(n^4+n^3+n^2)$ where $n$
is the dimension of the system. The algorithm is a generalisation of an algorithm presented in~\citet{AW04}.

For time-invariant systems, i.e., $A_k \equiv A, b_k \equiv b$ and $c_k \equiv c$, the condition $c \adj(A) b \ne 0$ is both necessary and sufficient
for memoryless linear output feedback nullification of a controllable and observable system, as shown in~\citet{AW04} (note that, for time-invariant
systems, the notions of controllable/observable and completely controllable/observable coincide).  For general time-varying systems, the situation is
more involved because it can be that $c_k \adj(A_k) b_k \ne 0$ for some $k$ values but not for all of them. The condition that this term vanishes for
every $k$, stated in the above theorem, is sufficient but not necessary for nullification.

In order to present the second result about discrete-time systems, the notions of controller canonical form and algebraic equivalence of systems are
needed. These notions are given by the following definitions.

\newcommand{\al}[1]{\alpha_{k,#1}}

\BD The system~\eqref{equ:system} is in a \emph{controller canonical form} if, for every $k \in \Z$, the matrix $A_k$ and the vector $b_k$ are of the
form
\begin{equation}
\label{equ:canonical-form}
A_k = \begin{pmatrix}
  0      & 1      &        & 0        \cr
  \vdots &        & \ddots &          \cr
  0      & 0      &        & 1        \cr
  \al{1} & \al{2} & \cdots & \al{n}
\end{pmatrix} \quad \mbox{and} \quad b_k=\begin{pmatrix}
  0 \cr
  \vdots \cr
  0 \cr
  1
\end{pmatrix}
\end{equation} where $\al{1},\al{2},...,\al{n}$ are scalars. \ED

\BD \label{def:equivalence} Two systems \system and \tildesystem are considered algebraically equivalent if there exists a sequence $\{T_k\}_{k \in \Z}$
of invertible transformations such that $\tilde{A}_k = T_{k+1} A_k {T_k}^{-1}$, $\tilde{b}_k = T_{k+1} b_k$ and $\tilde{c}_k = c_k {T_k}^{-1}$ for every
$k \in \Z$.  \ED

Given algebraically equivalent systems \system and \tildesystem\!\!, the sequence $\{x_k\}_{k \in \Z}$ obeys equations~\eqref{equ:system} if and only if
the sequence $\{\tilde{x}_k\}_{k \in \Z}$ defined by $\tilde{x}_k = T_k x_k$ obeys the equations
\begin{equation} \label{equ:canonical-system}
\begin{aligned} \tilde{x}_{k+1} &= \tilde{A}_k \tilde{x}_k + \tilde{b}_k u_k \\
\tilde{y}_k &= \tilde{c}_k \tilde{x}_k.
\end{aligned} \end{equation}
Therefore, trajectories of similar systems are translations of each other under the time-varying change of coordinates given by the sequence $\{T_k\}_{k
\in \Z}$.

\medskip

The second theorem concerning discrete-time systems is as follows.

\BT \label{thm:canonical-form} The system~\eqref{equ:system} is completely controllable if and only if it is algebraically equivalent to a system in the
controller canonical form~\eqref{equ:canonical-form}.\ET

In other words, \REF{Theorem}{thm:canonical-form} says that every completely controllable system can be represented by an $n$th order scalar equation.
This is a generalisation of the well known fact that every controllable time-invariant system can be represented as an $n$th order scalar equation.
\citet{Sil66} proved that continuous-time controllable time-varying systems with smooth coefficients can be represented by an $n$th order differential
scalar equation (see \REF{Remark}{rem:cannonical-form-of-cont-systems} below for a discussion related to this result). \citet{Gaishun00} used system
theoretic approach to prove a result similar to \REF{Theorem}{thm:canonical-form}. We provide a direct control theoretic proof of the theorem and give
an explicit formula for the coefficients in the canonical form. This theorem plays a focal role in the proof that we give for
\REF{Theorem}{thm:nullification}.

\medskip

Now, we present the results concerning continuous-time systems.

Consider a finite dimensional scalar-input scalar-output controllable and observable continuous-time time-varying system with analytic coefficients
\begin{equation} \label{equ:ct-system}
\begin{aligned} \dot{x}(t) &= A(t) x(t) + b(t) u(t) \\
y(t) &= c(t) x(t)
\end{aligned}\end{equation}
where $A(t)$, $b(t)$ and $c(t)$ are time-varying $n \times n$ matrix, $n$ dimensional column vector and $n$ dimensional row vector, respectively. We
assume that all the coefficients in these vectors are real analytic in $t$.

The sampling control strategy that we consider allows a feedback from the output. This strategy is to sample the output at prescribed equidistributed
times and to hold the control constant in periods between samplings. Let $\delta$ be the length of the sampling interval. For the discrete set of times
where sampling occur, say $t_k=k\delta, k \in \Z$, we get a discrete-time system of the form
\begin{equation} \label{equ:sampled-system}
\begin{aligned}
x_{k+1} &= A_{\delta}(k)x_k + b_{\delta}(k)u(k) \\
   y_k &= c_\delta(k) x_k
\end{aligned}
\end{equation}
where
\begin{equation} \label{equ:sampled-system-data}
\begin{aligned}
A_{\delta}(k) &=\Phi(t_{k+1},t_k), \\
b_{\delta}(k) &=\int_{t_k}^{t_{k+1}}\Phi(t_{k+1},s)b(s)\;ds, \\
c_\delta(k) &= c(t_k)
\end{aligned}
\end{equation}
and $\Phi$ is the fundamental matrix solution associated to $A(t)$.

For continuous-time systems, the control objective we consider is $\delta$-sample linear output feedback nullification as given in the next definition.

\BD We say that the continuous-time system~\eqref{equ:ct-system} is $\delta$-sample uniformly nullifiable by memoryless linear output feedback if the
sampled-data system~\eqref{equ:sampled-system} is uniformly nullifiable by memoryless linear output feedback. \ED

The first result concerning continuous-time systems is formulated in the following theorem. The term `almost every' means all except a countable set and
the term `generic' means that the property is valid for an open and dense set of observables, with respect to the supremum norm.

\BT \label{thm:CTVS-nullification} If the system~\eqref{equ:ct-system} is controllable and observable then, for a generic $c(t)$ and almost every
sampling period $\delta>0$, it is $\delta$-sample uniformly nullifiable by memoryless linear output feedback. \ET

The second result about continuous-time systems is the following preservation of controllability theorem.

\BT \label{thm:perservation-of-cont} If the system~\eqref{equ:ct-system} is controllable then, for almost every sampling period $\delta>0$, the
sampled-data system~\eqref{equ:sampled-system} is completely controllable. \ET

Conditions for preservation of controllability for time-invariant linear systems are studied in~\citet{KHN63} where the well known Kalman-Ho-Narendra
condition is presented. For time-invariant nonlinear systems, the problem is addressed in~\citet{Son83}.

The rest of this paper is organised as follows. In \REF{Section}{sec:controller-canonical-form} we prove \REF{Theorem}{thm:canonical-form} and add some
related remarks. \REF{Section}{sec:pers-of-cont} contains analysis of sampled-data systems and proofs of \REF{Theorem}{thm:CTVS-nullification} and
\REF{Theorem}{thm:perservation-of-cont}. In \REF{Section}{sec:nullification} a proof of \REF{Theorem}{thm:nullification} is presented.

\section{Controller canonical form}
\label{sec:controller-canonical-form} %
In this section we prove \REF{Theorem}{thm:canonical-form}. The proof is constructive and includes explicit formulas for the controller canonical
representation of a given system and for the transformation that brings the system to that form.

\medskip

We begin with the definition of the controllability matrices.

\BD For the control system~\eqref{equ:system} and $k \in \Z$, the $k$th controllability matrix  is defined as
\begin{equation}
\label{equ:controllability-matrix} %
W_k = \left[b_k,~A_k b_{k-1},~A_kA_{k-1} b_{k-2},~ ...~,~A_k \cdots A_{k-n+2}b_{k-n+1} \right].
\end{equation}
\ED


Invertibility of the controllability matrices corresponds to complete controllability of the system as shown in the next claim.

\BCM \label{cla:W_inv} The system~\eqref{equ:system} is completely controllable if and only if all the controllability matrices $\{W_k\}_{k \in \Z}$ are nonsingular.\ECM %
\BPF %
By~\eqref{equ:system}, if $x_{k-n+1}=0$ then $x_{k+1} = W_k (u_k,...,u_{k-n+1})^T$. Therefore, it is possible to steer
the zero state to any other state if and only if $W_k$ is nonsingular.
\EPF %

In the following definition we give a formula for the coefficients $\alpha_{i,j}$ that appear in the controller canonical form
\eqref{equ:canonical-form}. These coefficients are computed in terms of the controllability matrix of the system. For the moment, we just define these
numbers. Later on, we will show that these are the numbers appearing in the last rows of the matrices $\tilde{A}_k$ of the algebraically equivalent
controller canonical form representation.

\BD \label{def:canonical-coefficients} For a completely controllable system~\eqref{equ:system} and every $k \in \Z$, the controller canonical form
coefficients $\alpha_{k-i,i}$, $i=1,...,n$ are given by \begin{equation} \label{equ:alpha} (\alpha_{k,1},~\alpha_{k-1,2},~...~,~\alpha_{k-n+1,n})^T =
{W_k}^{-1}A_kW_{k-1}(0,...,0,1)^T
\end{equation} where $W_k$ denotes the $k$th controllability matrix~\eqref{equ:controllability-matrix}. \ED

Now we present the main proposition that states that the above numbers are indeed the coefficients in the controller canonical form. The transformation
that brings the system to this form is specified in the proof.

\BP \label{pro:canonical-form} Every completely controllable system~\eqref{equ:system} is algebraically equivalent to a system \tildesystem where
\[\tilde{A}_k =
\begin{pmatrix}
  0      & 1      &        & 0        \cr
  \vdots &        & \ddots &          \cr
  0      & 0      &        & 1        \cr
  \al{1} & \al{2} & \cdots & \al{n}
\end{pmatrix}, \quad \tilde{b}_k=\begin{pmatrix}
  0 \cr
  \vdots \cr
  0 \cr
  1
\end{pmatrix}
\] and $\al{1},\al{2},...,\al{n}$ are the scalars given in \REF{Definition}{def:canonical-coefficients}. \EP
\BPF%
Consider the controller canonical form controllability matrices
\begin{equation}
\label{equ:canonical-controllability-matrix}
 \tilde{W}_k = \left [\tilde{b}_k,~\tilde{A}_k
\tilde{b}_{k-1},~\tilde{A}_k \tilde{A}_{k-1} \tilde{b}_{k-2},~...~,~\tilde{A}_k \cdots
\tilde{A}_{k-n+2}\tilde{b}_{k-n+1} \right ].
\end{equation}

We first show that these matrices are nonsingular.

Note that $\tilde{A}_k$ is a shift matrix, i.e., for every $\xi_1,...,\xi_n \in \R$, $\tilde{A}_k (\xi_1,...,\xi_n)^T=(\xi_2,...,\xi_n,\sigma)$ for some
$\sigma \in \R$. In particular, the columns of $\tilde{W}_k$ are shifts of the vector $(0,...,0,1)^T$. This means that $\tilde{W}_k$ is a skew lower
triangular matrix with $1$'s in the main skew diagonal. Therefore, for every $k \in \Z$, the canonical controllability matrix $\tilde{W}_k$ is
nonsingular.

Now, we are ready to specify the transformation that brings the system into the controller canonical form:
\[ T_k = \tilde{W}_{k-1}^{ }W_{k-1}^{-1} \]
where $W_k$ and $\tilde{W}_k$ are the controllability matrices given in equations~\eqref{equ:controllability-matrix} and
\eqref{equ:canonical-controllability-matrix} respectively.  From the invertibility of $W_k$ and $\tilde W_k$ we get that $T_k$ is also invertible.

To complete the proof, we need to show that $T_{k+1}b_k = \tilde{b}_k$ and $T_{k+1}A_k{T_k}^{-1} = \tilde{A}_k$. The first part is easily verified:
\[T_{k+1}b_k=\tilde{W}_k{W_k}^{-1}b_k=\tilde{W}_k(1,0,...,0)^T=\tilde{b}_k.\]

By the definition of $W_k$ and $\tilde{W}_k$, the first $n\!-\!1$ columns of $\tilde{A}_k \tilde{W}_{k-1}$ are the same as the last $n\!-\!1$ columns of
$\tilde{W}_k$ and the first $n\!-\!1$ columns of $A_k W_{k-1}$ are the last $n\!-\!1$ columns of $W_k$. Therefore the first $n\!-\!1$ columns of
${W_k}^{-1}A_kW_{k-1}$ and ${\tilde{W}_k}^{-1}\tilde{A}_k\tilde{W}_{k-1}$ coincide. The last columns of these matrices are equal by equation
\eqref{equ:alpha} (since the last column of ${\tilde{W}_k}^{-1}\tilde{A}_k\tilde{W}_{k-1}$ is
$(\alpha_{k,1},~\alpha_{k-1,2},~...~,~\alpha_{k-n+1,n})^T$). We get that ${W_k}^{-1}A_kW_{k-1} = {\tilde{W}_k}^{-1}\tilde{A}_k\tilde{W}_{k-1}$ and
therefore
\[T_{k+1}A_k{T_k}^{-1}=\tilde{W}_k{W_k}^{-1}A_k\tilde{W}_{k-1}{W_{k-1}}^{-1}=\tilde{A}_k. \qedhere \] \EPF

The proof of \REF{Theorem}{thm:canonical-form} is given as a corollary of \REF{Claim}{cla:W_inv} and \REF{Proposition}{pro:canonical-form}.

\BPF[Proof of \REF{Theorem}{thm:canonical-form}] By the preceding proposition, every completely controllable system has a controller canonical form
representation. The converse follows from \REF{Claim}{cla:W_inv} and from the proof of \REF{Proposition}{pro:canonical-form} where it is shown that all
the controllability matrices of a system in a controller canonical form are nonsingular. \EPF

We conclude this section with two remarks. The first remark regards the relation of the results presented in this section with a result due to L. M.
Silverman, about the canonical form for continuous-time systems. And the second remark is about the use of negative time indices.

\BRM \label{rem:cannonical-form-of-cont-systems} The result stated as \REF{Theorem}{thm:canonical-form} above is a discrete-time analogue of a theorem
reported by~\citet{Sil66}.

For continuous-time time-varying linear systems with smooth coefficients of the form~\eqref{equ:ct-system}, the notions of complete controllability and
algebraic equivalence are as follows.

A continuous-time time-varying linear systems with smooth coefficients is said to be completely controllable if the controllability matrix
$W(t)=[p_0(t),p_1(t),...,p_{n-1}(t)]$ has full rank everywhere, where $p_{k+1}(t)=A(t)p_k(t)+\dot{p}_k(t)$ and $p_0(t)=b(t)$.

Two continuous-time time-varying linear systems $(A(t),b(t),c(t))_{t \in \R}$ and $(\tilde{A}(t),\tilde{b}(t),\tilde{c}(t))_{t \in \R}$ are said to be
algebraically equivalent if there exists a nonsingular time-varying matrix $T(t)$ with continuous derivative such that
$\tilde{A}(t)=(T(t)A(t)-\dot{T}(t)){T(t)}^{-1}$, $\tilde{b}(t)=T(t)b(t)$ and $\tilde{c}(t)=c(t)T(t)^{-1}$. The fundamental matrix solution of the
algebraic equivalent system is given by $\tilde\Phi(t,s)=T(t)\Phi(t,s)T(s)^{-1}$ where $\Phi(t,s)$ is the fundamental matrix solution of the original
system.

Note that algebraic equivalence and sampling commute: consider two algebraic equivalent smoothly varying systems $(A(t),b(t),c(t))_{t \in \R}$ and
$(\tilde{A}(t),\tilde{b}(t),\tilde{c}(t))_{t \in \R}$. Then, for every $\delta>0$, the sampled-data systems $(A_\delta(k),b_\delta(k),c_\delta(k))_{k
\in \Z}$ and $(\tilde{A}_\delta(k),\tilde{b}_\delta(k),\tilde{c}_\delta(k))_{k \in \Z}$, given by equations~\eqref{equ:sampled-system-data} for each
system respectively, are also algebraically equivalent (in the sense of \REF{Definition}{def:equivalence}). More specifically, if the equivalence
transformation is $T(t)$, then $\tilde{A}_\delta(k)=T(t_{k+1})A_\delta(k)T(t_k)^{-1}$, $\tilde{b}_\delta(k)=T(t_{k+1})b_\delta(k)$ and
$\tilde{c}_\delta(k)=c_\delta(k)T(t_k)^{-1}$.

\citet{Sil66} proved that a continuous-time time-varying linear system with smooth coefficients is algebraically equivalent to a system in the
controller canonical form if and only if it is completely controllable.

The resemblance of this  result to \REF{Theorem}{thm:canonical-form} is apparent. The proof in~\citet{Sil66} is also closely related to the proof of
\REF{Proposition}{pro:canonical-form} given here. In both proofs, the matrices $T(t)=\tilde{W}(t)^{}W(t)^{-1}$ are shown to be nonsingular under the
controllability assumption and then used to transform the system to its controller canonical form (where $\tilde{W}(t)$ is the controllability matrix of
the system in the controller canonical form).

This resemblance suggests that complete ($n$-step) controllability is a discrete analogue of complete controllability in the continuous-time case (thus
the name). Generally speaking, both notions of complete controllability relate to the ability to bring any initial state to any final state within every
interval.
 \ERM

\BRM In this paper we considered systems defined on the doubly infinite time domain $\Z$. When the time domain is not doubly infinite, complete
controllability and complete observability are assumed only for positive times (replace $\Z$ with $\N$ in both definitions). Note that , in that case,
complete controllability cannot guarantee controller canonical form at all positive times because, for example, it does not say anything about $A_0$.

However, except for a finite prefix, controller canonical form exists. From the proof of \REF{Proposition}{pro:canonical-form}, one can see that the
controller canonical form of $\tilde{A_n},\tilde{A}_{n+1},...$ is based on the invertibility of $W_0,W_1,...$ which corresponds to complete
controllability at positive times. Therefore, if we assume that any state at all nonnegative times can be steered to any other state in $n$ steps, there
exists a sequence of invertible transformations $\{ T_k \}_{k=n}^\infty$ such that $\tilde{A}_k=T^{}_{k+1}A^{}_kT^{-1}_k$ and $\tilde{b}_k=T_{k+1}b_k$
are in a controller canonical form for every $k \geq n$.

It is also possible to state a condition that guarantees transformation to a controller canonical form of all the matrices: if the matrices \[
[A_0,~b_0],~[A_1A_0,~A_1,~b_0],~...~,~[A_n \cdots A_0,~A_n \cdots A_1b_0,~...~,~A_nb_{n-1},~b_n] \] have full rank and any state at all nonnegative
times can be steered to any other state in $n$ steps, then there exists a sequence of invertible transformations $\{T_k\}_{k \in \N}$ such that
$\tilde{A}_k=T^{}_{k+1}A^{}_kT^{-1}_k$ and $\tilde{b}_k=T_{k+1}b_k$ are in a controller canonical form for all $k \in \N$. This claim is true since,
under the above condition, it is possible to extend the system backwards by adding matrices $A_{-n},...,A_{-1}$ and vectors $b_{-n},...,b_{-1}$ such
that any state at a time instance $k \geq -n$ can be steered to any other state at time $k+n$. \ERM

\section{Controllability under sampling} \label{sec:pers-of-cont}
The proof of \REF{Theorem}{thm:perservation-of-cont} is presented in two steps: first, a statement about analytic curves is given as
\REF{Proposition}{pro:pers-of-cont} and proved using some intermediate claims. Then the proof of the theorem is derived as a corollary of that
proposition.

\BP \label{pro:pers-of-cont} Let $\psi:\R \to \R^n$ be a real analytic curve. Assume that there is an uncountable 
set $\Delta \subseteq \R_{>0}$ such that, for every $\delta \in \Delta$, there exists $k_\delta \in \Z$ such that the vectors
\[\int_{k_\delta \delta}^{(k_\delta+1)\delta}\!\!\!\!\!{\psi(t)\;dt}, \quad \int_{(k_\delta+1)\delta}^{(k_\delta+2)\delta}\!\!\!\!\!\psi(t)\;dt, \quad ... \,, \quad \int_{(k_\delta+n-1) \delta}^{(k_\delta+n)\delta}\!\!\!\!\!\psi(t)\;dt\]
are linearly dependent. Then there exists a proper linear subspace $V \subseteq \R^n$ such that $\psi \big(\R\big) \subseteq V$. \EP

\BCM \label{cla:f_is_zero} Let $\psi:\R \to \R^n$ be a curve satisfying the conditions of \REF{Proposition}{pro:pers-of-cont}. Then there exists $k \in
\Z$ such that the function
\[f(\delta) = \det{\left[\int_{k \delta}^{(k+1)\delta}\!\!\!\!\!{\psi(t)\;dt}, \quad \int_{(k+1)\delta}^{(k+2)\delta}\!\!\!\!\!\psi(t)\;dt, \quad ...
\,, \quad \int_{(k+n-1) \delta}^{(k+n)\delta}\!\!\!\!\!\psi(t)\;dt \right]}\]
 vanishes for every $\delta >0$. \ECM
\BPF %
We have a map $\delta \mapsto k_\delta$ from an uncountable set to a countable set. By the Pigeonhole Principle, there must be $k \in \Z$ whose preimage
is uncountable. For this $k$, the real analytic function $f(\delta)$ has an uncountable zero set, therefore it is identically zero.
\EPF %

\BCM \label{cla:nth-der-of-f} Let $\psi:\R \to \R^n$ be a real analytic curve. The $m$th derivative at zero of the function defined in the preceding
claim is given by \[\sum_{\substack{m_1+\cdots+m_n=m \\ 0< m_1 < \cdots < m_n }} {m \choose m_1,...,m_n} \det\bigl(c(i,m_j)\bigr)_{i,j=1}^n
\det\left[\psi^{(m_1-1)}(0),...,\psi^{(m_n-1)}(0)\right]\] where $c(i,l)=(k+i)^l-(k+i-1)^l$ and $\psi^{(i)}(0)$ denotes the $i$th derivative of $\psi$
at zero. \ECM
\BPF %
Recall that the $m$th derivative of the determinant of a time-varying matrix, $M\colon \R \to \R^{n \times n}$ is given by
\begin{equation}\label{eqn:nth-der-of-det}
\sum_{m_1+\cdots+m_n=m}{m \choose m_1,...,m_n} \det \left[\frac{d^{m_1}}{dt^{m_1}} M_1,\frac{d^{m_2}}{dt^{m_2}}M_2,...,
\frac{d^{m_n}}{dt^{m_n}}M_n\right]
\end{equation}
where $M_i$ is the $i$th column of $M(t)$.

Consider the function $g_i(\delta)=\int_{(k+i-1)\delta}^{(k+i)\delta} \psi(t)\;dt$ which gives the $i$th column of the matrix inside the determinant in
$f(\delta)$ as a function of $\delta$. It is easy to verify that
\[ \frac{d^m}{d\delta^m}g_i(0) = \left\{%
\begin{array}{ll}
    c(i,m) \psi^{(m-1)}(0), & m>0; \\
    0,                 & m=0. \\
\end{array}%
\right. \] %
In particular, by formula~\eqref{eqn:nth-der-of-det}, the $m$th derivative of $f$ is
\[\sum_{\substack{m_1+\cdots+m_n=m \\ m_1 m_2 \cdots m_n \ne 0}} {m \choose m_1,...,m_n} \det\left[c(1,m_1)\psi^{(m_1-1)}(0),...,c(d,m_n)\psi^{(m_n-1)}(0)\right].\] %
Factoring out scalars from the columns yields
\[\sum_{\substack{m_1+\cdots+m_n=m \\ m_1 m_2 \cdots m_n \ne 0}} {m \choose m_1,...,m_n} \prod_{i=1}^n c(i,m_i)  \det\left[\psi^{(m_1-1)}(0),...,\psi^{(m_n-1)}(0)\right].\] %
Collecting together terms that corresponds to permutations of the same partition gives
\[
\sum_{\substack{m_1+\cdots+m_n=m \\ 0< m_1 < \cdots < m_n }} {m \choose m_1,...,m_n}  \sum_{\pi \in S_n} \prod_{i=1}^n c(i,m_{\pi(i)}) \sgn(\pi) \det\left[\psi^{(m_1-1)}(0),...,\psi^{(m_n-1)}(0)\right].  %
\]
The claimed formula follows from the definition of the determinant.
\EPF %

We proceed with the proof of \REF{Proposition}{pro:pers-of-cont}. The proof scheme is to use the fact that $f$ is identically zero in order to prove
that the range of $\psi$ is confined within a proper linear subspace. We will do this by proving that the derivatives of $\psi$ at the origin are all in
a proper linear space using the fact that all the derivatives of $f$ are zero. The main tool for this proof scheme is provided by the following lemma.

\BL \label{Lem:rows-implies-all} Let $x_1,x_2,...$ be a sequence of $n$ dimensional vectors. Assume that there is a function $C: \N^n \to \R
\smallsetminus \{0\}$ such that
\begin{equation}\label{equ:assump}
\sum_{\substack{m_1+\cdots+m_n=m \\ 0< m_1 < \cdots < m_n }}C(m_1,...,m_n) \det[x_{m_1},x_{m_2},...,x_{m_n}] = 0
\end{equation} for every $m \in \N$. Then the sequence is contained in a proper linear subspace. \EL
\BPF %
We begin by introducing a linear order over the set of ordered tuples. We write $(m_1,m_2,...,m_n) \prec (\hat{m}_1,\hat{m}_2,...,\hat{m}_n)$ if
$\left(\sum_{i=1}^n m_i,\hat{m}_n,\hat{m}_{n-1},...,\hat{m}_1\right)$ precedes $\left(\sum_{i=1}^n\hat{m}_i,m_n,m_{n-1},...,m_1\right)$
lexicographically. Note that this is a well-founded order~\citep{Wei03}.

Let $(m_1,m_2,...,m_n)$ be the minimal (according to the above order) ordered tuple for which
\begin{equation}\label{equ:det-is-not-zero}
\det[x_{m_1},x_{m_2},...,x_{m_n}] \ne 0.
\end{equation}

Take another ordered tuple, $(\hat{m}_1,\hat{m}_2,...,\hat{m}_n)$, such that $\sum_{i=1}^n m_i = \sum_{i=1}^n \hat{m}_i$. Towards a contradiction to the
existence of a tuple satisfying equation~\eqref{equ:det-is-not-zero}, we will show that $\det[x_{\hat{m}_1},x_{\hat{m}_2},...,x_{\hat{m}_n}]=0$, i.e.,
that all the terms in equation~\eqref{equ:assump} vanish, except the one that corresponds to $(m_1,m_2,...,m_n)$.

Consider first the case where $(m_n,m_{n-1},...,m_1)$ precedes $(\hat{m}_n,\hat{m}_{n-1},...,\hat{m}_1)$ lexicographically. Thus,
$(\hat{m}_n,\hat{m}_{n-1},...,\hat{m}_1)$ precedes $(m_n,m_{n-1},...,m_1)$ in our order. Since $(m_1,m_2,...,m_n)$ is
the first tuple for which the determinant is not zero, we have $\det[x_{\hat{m}_1},x_{\hat{m}_2},...,x_{\hat{m}_n}]=0$. %

Assume that $(\hat{m}_n,\hat{m}_{n-1},...,\hat{m}_1)$ precedes $(m_n,m_{n-1},...,m_1)$ lexicographically. Let $i \in \{0,1,...,n-1\}$ be the first index
such that $m_{n-i} \ne \hat{m}_{n-i}$ (more specifically, $m_{n-i} > \hat{m}_{n-i}$). Then, $\hat{m}_{n-j} = m_{n-j}$ for $j=0,...,i-1$ and
$\hat{m}_{n-j} < m_{n-i}$ for $j=i,...,n-1$. Therefore,
\[\det[x_{m_1},...,x_{m_{n-i-1}},x_{\hat{m}_j},x_{m_{n-i+1}},...,x_{m_n}] = 0\] for all $j=1,2,...,n$ (some because of repeated columns and the others
because $\hat{m}_j- m_{n-i}+\sum_{k=1}^n m_k < \sum_{k=1}^n m_k$). Since the vectors $x_{m_1},x_{m_2},...,x_{m_n}$ are linearly independent
(by~\eqref{equ:det-is-not-zero}), we get that
\[\{x_{\hat{m}_1},x_{\hat{m}_2},...,x_{\hat{m}_n}\} \subset \spn \left(\{x_{m_1},x_{m_2},...,x_{m_n}\} \smallsetminus \{x_{m_{n-i}}\}\right).\]
In particular, $\det[x_{\hat{m}_1},x_{\hat{m}_2},...,x_{\hat{m}_n}]=0$.

We get that, all the terms in equation~\eqref{equ:assump}, except the term that corresponds to $(m_1,m_2,...,m_n)$, vanish. This yields a contradiction
with inequality~\eqref{equ:det-is-not-zero}.
\EPF %

Note the similarity of equation~\eqref{equ:assump} and the expression for the $m$th derivative of $f$ given in~\eqref{cla:nth-der-of-f}. To apply
\REF{Lemma}{Lem:rows-implies-all}, we need to verify that the coefficients are not zero. In the following claim we show that they are all positive.

\BCM \label{cla:positive-coeff} For every $k>0$ and integers $0<m_1<m_2<\cdots<m_n$, the matrix \[ M = \big(c(i,m_j)\big)_{i,j=1}^n = \big(
(k+i)^{m_j}-(k+i-1)^{m_j} \big)_{i,j=1}^n
\]
has a positive determinant. \ECM
\BPF %

Define
\[f(k\,;m_1,m_2,...,m_n)= \det M.\]
Note that this functions vanishes when there is a repeating parameter, i.e., $m_i=m_j$ for some $i \ne j \in \{0,1,..,n\}$ where $m_0=0$.

By adding the rows of $M$ and deleting telescopic terms, it is easy to see that
\[ {\small \begin{pmatrix}
1      & 0      & \cdots & 0 \cr %
1      & 1      & \ddots & \vdots \cr %
\vdots & \vdots & \ddots & 0 \cr %
1      & 1      & \cdots & 1 \cr %
\end{pmatrix}}  M  = \big(
(k+i)^{m_j}-k^{m_j} \big)_{i,j=1}^n\,. \]%
When $k\!=\!0$, this is a generalised Vandermonde matrix \citep{Wei03}. Since the determinant of a generalised Vandermonde is positive, we get that
$f(0\,;m_1,m_2,...,m_n)>0$.

By formula~\eqref{eqn:nth-der-of-det}, the derivative of $f(k\,;m_1,m_2,...,m_n)$ with respect to $k$ is
\begin{equation}
\label{equ:der-of-f} \frac{d}{dk}f(k\,;m_1,m_2,...,m_n)=\sum_{i=1}^n m_i f(k\,;m_1,...,m_i-1,...,m_n).
\end{equation}

We now prove, by induction on $m = m_1+m_2+\cdots+m_n$, that $f(k\,;m_1,m_2,...,m_n)>0$ for all $k>0$ and integers $0<m_1<m_2<\cdots<m_n$.

By~\eqref{equ:der-of-f}, $\frac{d}{dk}f(k\,;1,2,...,n)=0$ for every $k$ (because all the terms have repeating parameter: $m_{i-1}=m_i-1$). Therefore,
$f(k\,;1,2,...,n)=f(0\,;1,2,...,n)>0$. This establishes the case $m=n(n+1)/2$, which is the base of the induction.

If $m_1+m_2+\dots+m_n>n(n+1)/2$ then each of the summands on the right hand side of~\eqref{equ:der-of-f} is nonnegative, by the induction hypothesis, so
the derivative $\frac{d}{dk}f(k\,;m_1,m_2,...,m_n)$ is positive. In particular, for every $k>0$, $f(k\,;m_1,m_2,...,m_n) \geq f(0\,;m_1,m_2,...,m_n)>0$.
\EPF %

The proof of the proposition follows from the preceding claims.

\BPF[Proof of \REF{Proposition}{pro:pers-of-cont}] By \REF{Claim}{cla:f_is_zero} and \REF{Claim}{cla:nth-der-of-f}, the sequence of the derivatives at
zero: $\psi(0),\psi^{(1)}(0),...$ satisfies the conditions of \REF{Lemma}{Lem:rows-implies-all} where the constants are given by $C(m_1,m_2,...,m_n) ={m
\choose m_1,...,m_n} \det\big(c(i,m_j)\big)_{i,j=1}^n$. By \REF{Claim}{cla:positive-coeff}, these coefficients are all positive hence all the
derivatives of $\psi$ at the origin lie in a proper linear subspace. Since $\psi$ is analytic, its whole image is contained in that subspace. \EPF

The proof of the theorem follows as a corollary of the preceding proposition.

\BPF[Proof of \REF{Theorem}{thm:perservation-of-cont}] Consider the curve $\psi(t)=\Phi(0,t)b(t)$, where $\Phi$ is the fundamental matrix solution
associated to $A(t)$. Let $\Delta$ be the set of $\delta$ values for which $\Sigma_{[\delta]}$ is not completely controllable. By
\REF{Claim}{cla:W_inv}, for every $\delta \in \Delta$ there exists $k_\delta \in \Z$ such that the $k_\delta$'s controllability
matrix~\eqref{equ:controllability-matrix} of the sampled data system~\eqref{equ:sampled-system}, namely, the matrix
 \[ \Phi\big((k_\delta+1)\delta,0\big) \left[ \int_{k_\delta \delta}^{(k_\delta+1)\delta}\!\!\!\!\!{\psi(t)\;dt}, \quad
\int_{(k_\delta-1)\delta}^{k_\delta\delta}\!\!\!\!\!\psi(t)\;dt, \quad ... \,, \quad \int_{(k_\delta-n+1)
\delta}^{(k_\delta-n+2)\delta}\!\!\!\!\!\psi(t)\;dt \right]   \] is singular. If $\Delta$ is not countable, the conditions of
\REF{Proposition}{pro:pers-of-cont} are met, so the image of $\psi$ is contained in a proper subspace. In particular, $\Sigma$ is not controllable
\citep[page 109, Theorem 5]{Son98}. \EPF

For the nullification algorithm presented in this paper, in addition to complete controllability and complete observability of the discrete-time system,
we need that $c_k \adj(A_k) b_k \ne 0$ for every $k \in \Z$ (see \REF{Theorem}{thm:nullification}). The following example shows that it may be that, for
all sampling periods, this condition is not satisfied; even if the continuous-time system is controllable and observable with analytic coefficients.

\BEX An example of a controllable and observable system with analytic coefficients such that there exists $k \in \Z$ for which $c_\delta(k)
\adj(A_\delta(k)) b_\delta(k)=0$ for every sampling period $\delta>0$. Choose an arbitrary $k$. Take $A(t) \equiv
0,b(t)=(1\!-\!n,2t,3t^2,...,nt^{n-1})^T$ and $c(t)=(\lambda_1,\lambda_2,...,\lambda_n)$ where $\lambda_i=-(t/k)^{n-i}((k+1)^i-k^i)^{-1}$. \EEX

From a practical point of view, imagine that we have a continuous-time system and choose some sampling period for which the sampled-data system is
controllable. We know, from the above example, that the system may not satisfy the sufficient condition for nullification. However, for any fixed $k \in
\Z$, almost any sequence of observation vectors allows to steer any initial state at time $k$ to the origin in a finite number of steps, as shown in the
following proposition.

\BP \label{pro:gen-obs-null} Consider the continuous-time system~\eqref{equ:ct-system}. If the sampled-data system~\eqref{equ:sampled-system} is
completely controllable then there exists $N \in \N$ such that for every $k \in \Z$ and almost every $c_k,...,c_{k+N} \in \R^{1 \times n}$, there are
scalars $F_k,...,F_{k+N-1}$ such that the sequence $x_k, x_{k+1},...,x_{k+N}$ resulting from the dynamics~\eqref{equ:dynamics} satisfies $x_{k+N}=0$.
\EP
\BPF%
Let $N=2(n^4+n^3+n^2)$. By \REF{Claim}{cla:W_inv}, if the sampled-data system is completely controllable then $b_\delta(k) \ne 0$ for every $k \in \Z$.
In particular, since the matrices $\adj(A_\delta(k))$ are nonsingular, we also have that $\adj(A_\delta(k))b_\delta(k) \ne 0$. Therefore, the set of
observables $c_k,...,c_{k+N} \in \R^{1 \times n}$ for which $c_i\adj(A_\delta(i))b_\delta(i) \ne 0$ for every $k \leq i \leq k+N$, is the a finite
product of complements of hyperplanes. In particular, if we intersect this set with the set of observables that yields a completely observable system we
get a set of measure one. By \REF{Theorem}{Thm:null-with-bound}, the systems in this set are nullifiable
\EPF%

The number $N$ in the above proposition is the same as in \REF{Definition}{def:nullification}. In particular, as in \REF{Theorem}{Thm:null-with-bound},
it is bounded by $2(n^4+n^3+n^2)$ where $n$ is the dimension of the system.

\medskip

We conclude this section by a proof of \REF{Theorem}{thm:CTVS-nullification}.

\BPF[Proof of \REF{Theorem}{thm:CTVS-nullification}]%
By \REF{Theorem}{thm:perservation-of-cont}, the sampled-data system~\eqref{equ:sampled-system} is completely controllable and completely observable for
almost any sampling period (using the duality principle). In that case, by \REF{Theorem}{thm:nullification}, \[\{ c(\cdot)\colon c(k
\delta)\adj(A_\delta(k))b_\delta(k) \ne 0 \mbox{ for every } k \in \Z\}\] is a subset of observables for which the system is uniformly nullifiable by
memoryless linear output feedback. By \REF{Claim}{cla:W_inv}, if the sampled-data system is completely controllable then $b_\delta(k) \ne 0$ for every
$k \in \Z$. Since the matrices $\adj(A_\delta(k))$ are nonsingular, we also have that $\adj(A_\delta(k))b_\delta(k) \ne 0$. In particular, the above set
of observables consists of the functions that avoid a sequence of $(n-1)$-dimensional hyperplanes on a discrete set of times. It is easy to verify that
such a set is open and dense in the uniform topology.
\EPF%

\section{Nullification by memoryless output feedback } \label{sec:nullification}
In this section we prove \REF{Theorem}{thm:nullification}. The proof is a generalisation of the proof of Theorem D. presented in~\citet{AW04} where
time-invariant systems are analysed.

We begin with a proposition that allows to consider only systems in a controller canonical form.

\BP \label{pro:nullification-invariant} If the system \system is algebraically equivalent to the system \tildesystem and \system is uniformly
nullifiable by memoryless linear output feedback then \tildesystem is also uniformly nullifiable by memoryless linear output feedback. \EP
\BPF %
Recall equations~\eqref{equ:canonical-system}. Because $\tilde{x}_k$ is defined as the image of $x_k$ under a bijective linear transformation, $x_k$ is
steered to the origin if and only $\tilde{x}_k$ is.
\EPF %

Because \REF{Theorem}{thm:nullification} is only about completely controllable systems and because every completely controllable system have a
controller canonical form representation (\REF{Theorem}{thm:canonical-form}), we will assume from now on that the system is given in a controller
canonical form.

To simplify notations, we will drop the tildes and write $A_k,b_k,c_k$ instead of $\tilde{A}_k,\tilde{b}_k,\tilde{c}_k$ respectively, keeping in mind
that the data is assumed to be in a controller canonical form~\eqref{equ:canonical-form}. We also assume that nullification begins in time zero. This
assumption imposes no loss of generality since it is always possible to shift time.

The first step towards a proof of \REF{Theorem}{thm:nullification} is the following proposition. There are three differences between this proposition
and the theorem. The first difference is that the theorem deals only with systems in a controller canonical form. The second difference is that in the
proposition the initial state is given, where in the theorem the same feedback must fit all initial states. The third difference is that we start
nullification at time zero and not at any time.

\BP \label{pro:specific-state-nullification} Consider a control system~\eqref{equ:system} represented in a controller canonical form
\eqref{equ:canonical-form} such that, for every $k \in \N$,  $c_k (1,0,...,0)^T \ne 0$ (namely, the first coordinate of $c_k$ is not zero). Then there
is a natural number $N \in \N$ such that for any initial state $x_0 \in \R^n$ there are coefficients $F_0,F_1,...,F_N$ such that the feedback $u_k = F_k
y_k$ achieves $x_N=0$. \EP

Towards a proof of \REF{Proposition}{pro:specific-state-nullification}, for a system satisfying the conditions of the proposition, consider the
following construction. The idea is to encode the next state relation of the system as an affine formula, unroll this to finite time, and to analyse the
resulting sequence.

\BCON \label{con:main} Starting with an arbitrary vector $x_0 = (\xi_{0,1},...,\xi_{0,n})^T$, the sequence $x_0,x_1,...$ is generated as follows:
\begin{itemize}
    \item If $c_0 x_0 = 0$ define $x_1 = A_0x_0$.
    \item If $c_0 x_0 \ne 0$ define $x_1 = (\xi_{0,2},...,\xi_{0,n},\delta_1)^T$ where $\delta_1$ is a variable whose value will be determined later.
\end{itemize}
Inductively, suppose that $x_k=(\xi_{k,1},...,\xi_{k,n})^T$ has been constructed.
\begin{itemize}
    \item If $c_k x_k = 0$ for any choice of numerical value of the variables $\{\delta_i \colon i \leq k\}$ define $x_{k+1} = A_{k}x_k$.
    \item Otherwise, introduce a new free variable $\delta_{k+1}$ and  define the next vector by $x_{k+1}=(\xi_{k,2},...,\xi_{k,n},\delta_{k+1})^T$.
\end{itemize}
 \ECON

Note that not all the variables in $\{\delta_i \colon i \in \N \}$ affect the coordinates of the vectors. The following notation is used to refer to the
variables that need to be assigned with a numerical value in order to make the trace concrete.

\BNOT \label{not:active} If $c_{k-1}x_{k-1} \ne 0$ for some numerical realisation of $\{\delta_i \colon i < k \}$ then the free variable $\delta_k$ is
called \emph{active}. A coordinate $i \in \{1,...,n\}$ of a vector $x_k$ is called active if the variable $\delta_{k-n+i}$ is active. \ENOT

The coordinates of the vectors $x_0,x_1,...$ introduced along the sequence presented in the above construction are affine formulas in the active
variables. Our proof scheme is to find assignment to these variables such that the last vector of the sequence is zero and all the other vectors are a
trace of the system under some feedback.

A focal object in the analysis is the sequence $d(0),d(1),...$ that counts the number of active coordinates in the vectors $x(0),x(1),...$ defined in
\REF{Construction}{con:main}. The number $d(k)$ have several interpretations as follows.

\BNOT Let $d(k)$ denote the number of active variables in $\{ \delta_i \colon k-n< i \leq k \}$. Equivalently, $d(k)$ is the number of active coordinates
of $x_k$. Another interpretation of $d(k)$ is the number of indices in $\{ i \colon k-n< i \leq k \}$ for which $c_{i-1}x_{i-1} \ne 0$ for some
assignment of the free variables. \ENOT

The properties given in the following three claims are the reason for our interest in the above sequence.

\BCM \label{cla:active-or-zero} For every $k \in \N$, $d(k) \leq d(k+n)$. Furthermore, if $d(k)=d(k+n)$ for every $k_0 \leq k < k_0+2n$ then the
coordinates of the vectors $x_{k_0},x_{k_0+1},..,x_{k_0+n}$ are either active or zero. \ECM
\BPF %
For a given $k \in \N$, consider the finite prefix $x_0,x_1,...,x_{k-n}$. The variables in this prefix are $\{ \delta_i \colon i \leq k-n \}$. For these
variables, fix a numerical realisation such that $c_{i-1}x_{i-1} \ne 0$ whenever $\delta_i$ is active. Note that, under such a realisation, the prefix
is concrete, i.e., all the entries are fixed numbers without free variables.

We claim that such a realisation exists: consider the Euclidian space of numerical realisations of the free variables introduced in the first $k$ steps
of \REF{Construction}{con:main}. Denote this space by $\R^{f(k)}$ where $f(k)$ is the number of free variables introduced until the $k$th step. For
every $i<k$, the set of realisations such that $c_ix_i \ne 0$ is the complement of an affine subspace in $\R^{f(k)}$. Therefore, if not empty, it must
be an open dense set. By definition, $\delta_{i+1}$ is active only if this set is not empty. In a finite prefix, it is possible to find a realisation
that satisfies $c_ix_i \ne 0$ whenever $\delta_{i+1}$ is active because the intersection of open dense sets is not empty.

Consider also the extended prefix, $x_0,x_1,...,x_k,x_{k+1},...,x_{k+n}$, under the same realisation. The first $k-n$ vectors are fixed whence the last
$2n$ may contain active variables and affine functions of active variables.

Now, the realisation is extended by fixing also the active variables introduced in $x_{k+1},...,x_{k+n}$. If $i>k$ and $\delta_i$ is active, set
$\delta_i = a_{i-1}x_{i-1}$ (where $a_{i-1}$ is the last row of $A_{i-1}$). We are left with only the $d(k)$ active coordinates of $x_k$ as free
variables. Denote the linear space of the realisations of these free variables by $\R^{d(k)}$.

Note that for every realisation in $\R^{d(k)}$, the vectors $x_k,x_{k+1},...,x_{k+n}$ are a trace of the autonomous system $x_{i+1}=A_ix_i$. To see
this, let $i \in \{k,k+1,...,k+n-1\}$. If $c_ix_i=0$ then $x_{i+1}=A_ix_i$ by \REF{Construction}{con:main}. Otherwise, the last entry of $x_{i+1} $ is
equal to the last entry of $A_ix_i$ by the extension of the realisation described above. The other entries of $x_{i+1}$ must agree with the
corresponding entries of $A_ix_i$ because of the shift structure of \REF{Construction}{con:main} and the controller canonical form of $A_i$.

Let $L_k$ be the mapping which assigns to an element in $\R^{d(k)}$ the string $\{c_kx_k,...,c_{k+n-1}x_{k+n-1}\}$ of observations. Only $d(k+n)$ of
these observations are not identically zero (by the last part of \REF{Notation}{not:active}). Thus, $L_k$ is considered as a mapping from $\R^{d(k)}$ to
the linear space $\R^{d(k+n)}$ of those $i$'s where $c_{i-1}x_{i-1}$ is not guaranteed to vanish.

If the system is observable then $L_k$ is one to one. For, if the mapping $L_k$ is not one to one, there are two realisations of the free variables in
$\R^{d(k)}$ which give rise to two distinct dynamics of length $n+1$ of the autonomous system $x_{i+1}=A_ix_i$ with the same observations.

This proves the first part of \REF{Claim}{cla:active-or-zero} because an affine mapping cannot be one to one if the dimension of the range is smaller
than the dimension of its domain.

\medskip

Towards a proof of the second part of the proposition, note that if $d(k)=d(k+n)$ then the mapping $L_k$ is one to one and onto because the dimension of
its range equals the dimension of the domain.

Assume that $d(k)=d(k+n)$. Because $L_k$ is onto, the zeroes observation is included in its range. Observability implies that an all zeroes observation
can only come from a null initial vector. Because $L_k$ is one to one, all the entries in $x_k$ must vanish if the free variables in that vector are set
to zero. The conclusion is that the entries of $x_k$ are linear (not only affine) functions of the active coordinates of $x_k$.


If $d(i)=d(i+n)$ for every $j \leq i < j+n$ then the coordinates of $x_j$ are either active or zero. This is true since, because of the shift structure,
every coordinate of $x_j$ becomes first in some vector $x_i$, $i \leq j < i+n$. The first entry cannot depend on variables introduced later in the
process so the only possible linear functions are constant zero or an active variable.

The second part of the proof of \REF{Claim}{cla:active-or-zero} follows by applying the above claim for $j=k,k\!+\!1,...,k\!+\!n$.
\EPF %

The following lemma provides a tool to extract information about the sequence $d(0),d(1),...$ from analysis of the sequences $d(i),d(i+n),d(i+2n),...$,
for $i=0,1,...,n-1$.


\BL \label{lem:eq-in-n-step-implies-dec} If $d(k + n) = d(k)$ then $d(k+n+1) \geq d(k+n)$. \EL
\BPF%
If the claim is false then $d(k+n) > d(k+n+1)$, i.e., the number of active coordinates decreases at step $k+n$. By \REF{Claim}{cla:active-or-zero},
$d(k+n+1) \geq d(k+1)$ hence the condition $d(k+n)=d(k)$ implies that $d(k) > d(k+1)$, i.e., the number of active coordinates decreases also at step
$k$. The shift structure implies that the number of active variables decreases at step $i$ only if the first coordinate of $x_i$ is active and the last
coordinate of $x_{i+1}$ is not active. For $i=k$ we get that the last coordinate of $x_{k+1}$ is not active, and for $i=k+n$ we get that the first
coordinate of $x_{k+n}$ is active. Since the last coordinate of $x_{k+1}$ is the first coordinate of $x_{k+n}$ we have a contradiction.
\EPF%

Towards the application of the second part of \REF{Claim}{cla:active-or-zero}, the following two claims give properties of the sequence $d(k)$.

\BCM \label{cla:existence-of-k0} There exists $k_0 \leq n^3+n^2$ such that $d(k)=d(k+n)$ for every $k_0 \leq k \leq k_0+n$. \ECM
\BPF %
Define the set $K=\{k\colon d(k) \ne d(k+n \}$. For  $i\!=\!0,...,n-1$ consider the sequence $d(i),d(i+n),d(i+2n),...$ which is nondecreasing
(\REF{Claim}{cla:active-or-zero}) and bounded by $n$. There are $n$ such sequences, each sequence increases at most $n$ times, so the size of $K$ is at
most $n^2$.

Assume that there exists no $k_0 < n^3+n^2$ such that $d(k)=d(k+n)$ for every $k_0 \leq k \leq k_0+n$. In particular, for all the intervals $I_j=j(n+1)
+ \{0,1,...,n\}$, $j=0,...,n^2-1$; the intersections $K \cap I_j$ are not empty. Therefor, the number of elements in $K \cap \{0,1,...,n^3+n^2-1\}$ is
at least $n^2$. In that case, since the size of $K$ is at most $n^2$, $K$ is bounded by $n^3+n^2$ so the claim is true for $k_0=n^3+n^2$.
\EPF %

\BCM \label{cla:existence-of-k0-2} There exists $k_0 \leq n^3+n^2$ such that $d(k)=d(k+1)$ for every $k_0 \leq k < k_0+2n$. \ECM
\BPF %
By \REF{Claim}{cla:existence-of-k0}, there exists $k_0 \leq n^3+n^2$ such that $d(k)=d(k+n)$ for every $k_0 \leq k \leq k_0+n$. By
\REF{Lemma}{lem:eq-in-n-step-implies-dec}, $d(k+n+1) \geq d(k+n)$ for every $k_0 \leq k \leq k_0+n$. Since $d(k_0+n)=d(k_0+2n)$, these inequalities
collapse to the equality $d(k_0+n)=d(k_0+n+1)=\cdots=d(k_0+2n)$. Using the equalities given by \REF{Claim}{cla:existence-of-k0} again, we get that
$d(k_0)=d(k_0+1)=\cdots=d(k_0+2n)$.
\EPF %

Using the second part of \REF{Claim}{cla:active-or-zero}, we now translate the property of $d(k)$ revealed in the previous claim, to properties of the
vectors $x(k)$ introduced in \REF{Contraction}{con:main}.

\BCM \label{cla:def_of_k_0} There exists $k_0 \leq n^3+n^2$ such that for every $k_0 \leq  k \leq k_0+n$:
\begin{enumerate}
    \item The entries of $x_k$ are either zero or active variables (no nonzero constants or affine functions).
    \item If a new variable is introduced in $x_{k+1}$ ($\delta_{k+1}$ is active) then the first entry of $x_k$ is a free variable ($\delta_{k-n+1}$ is also active).
\end{enumerate}
\ECM
\BPF%
By \REF{Claim}{cla:existence-of-k0-2} there exists $k_0 \leq n^3+n^2$ from which $d(k)$ is constant for $2n$ consecutive indices. By the second part of
\REF{Claim}{cla:active-or-zero}, if the sequence $d(k)$ is constant for $2n$ consecutive indices then the entries in the vectors in these indices are
either zero or free variables.

If $d(k)=d(k+1)$ then the number of free variables in $x_k$ equals the number of free variables in $x_{k+1}$. In particular, a new variable is
introduced in $x_{k+1}$ only if there is a free variable in $x_k$ which is not in $x_{k+1}$. Because of the shift structure, this can only happen if the
first entry of $x_k$ is a free variable.
\EPF%

{ \samepage With reference to the number $k_0$, identified in the preceding claim, we fix a new realisation of the free variables such that:
\begin{enumerate}
    \item For every $i < k_0 +n$ such that $\delta_{i+1}$ is active, $c_ix_i \ne 0$.
    \item For every $i>k_0$, $\delta_i=0$.
\end{enumerate} }

\BCM Such a realisation exists. \ECM
\BPF%
As in the proof of \REF{Claim}{cla:active-or-zero}, the set of realisations in $\R^{f(k_0-n)}$ satisfying $c_ix_i \ne 0$ for every $i \leq k_0$ for
which $\delta_{i+1}$ is active, is an intersection of open dense set and therefore not empty.

For $i > k_0$: by the first part of \REF{Claim}{cla:def_of_k_0}, the term $c_ix_i$ is linear in the active variables. Moreover, the coefficient of
$\delta_{i-n}$ in that term is the first coordinate of $c_i$ which is not zero by assumption.

By the second part of \REF{Claim}{cla:def_of_k_0}, if $\delta_i$ is active then $\delta_{i-n}$ is also active. Thus, for every $i>k_0$ for which
$\delta_i$ is active we can use the freedom in the variables $\delta_{k_0-n},...,\delta_{k_0}$ to make $c_ix_i$ not zero.

The second condition does not contradict the first one because the variables in $\{ \delta_i \colon i > k_0 \}$ have no affect on the numbers $\{ c_i x_i
\colon 0 \leq i < k_0 + n \}$ .
\EPF%

Since the last $n$ steps are shifts with $0$ entering in the last coordinate, it is clear that the vector $x_{k_0+n}$ is
zero. To finish the proof of \REF{Proposition}{pro:specific-state-nullification} we need to show that the vectors
$x_0,...,x_{k_0+n}$ are generated as the trace of the system under a controller of the form $u_k=F_k y_k$.

\BPF[Proof of \REF{Proposition}{pro:specific-state-nullification}] Using the above realisation,  define the feedback coefficients
\[
    F_k = \begin{cases}
        0,   &   \text{if $c_k x_k = 0$}; \\
        \frac{\delta_{k+1} + a_k x_k}{c_k x_k}, &   \text{otherwise}.
    \end{cases}
\]
where $a_k$ is the last row of $A_k$. By \REF{Construction}{con:main}, the feedback $u_k = F_k y_k$ generates
$x_0,x_1,...,x_{k_0+n}$ as a trace. In particular, it steers $x_0$ to the origin in finite time. \EPF

In \REF{Proposition}{pro:specific-state-nullification} we only assert that given a state $x_0 \in \R^n$, there exists a feedback that steers $x_0$ to
the origin. To prove \REF{Theorem}{thm:nullification} we need to swap the quantifiers, i.e., to show that there is a feedback that steers all initial
states to the origin. It turns out that these properties are equivalent, as shown in the following proposition.

\BP \label{pro:enough_to_nullify_a_given_state}  Given a control system~\eqref{equ:system}. Suppose that there is $N \in \N$ such that for every $\xi
\in \R^n$ and every $k \in \N$ there are $F_k,...,F_{k+N} \in \R$ such that the initial state $x_k=\xi$ with the controller $u_i=F_iy_i$ give
$x_{k+N}=0$. Then the system is output feedback nullifiable. \EP
\BPF %
Let $W_i  = (A_i+F_i b_i c_{i-1})$. By equations~\eqref{equ:system}, given $u_i=F_iy_i$, we have $x_{i+1}=W_i x_{i}$. Thus, the matrix $W_{k+N} \cdots
W_k$ maps the $k$th state to the $(k\!+\!N)$th state. This mapping is parameterised by $F_k,...,F_{k+N}$. We are given that for every $k \in \N$ and $x
\in \R^n$ there are $F_k,...,F_{k+N} \in \R$ such that $W_{k+N} \cdots W_kx=0$. Call that matrix $H(x,k)$, i.e., $H(x,k) x = 0$.

Let $v_1,...,v_n$ be a basis. Let $M_1 = H(v_1,k)$ and $M_{i+1} = H(M_i \cdots M_1 v_{i+1},k+im)$. The product $M=M_n \cdots M_1$ satisfies, $Mv_i =0$
for every $i=1,...,n$. Since $v_1,...,v_n$ is a basis we get that $M$ is the zero matrix.

The matrix $M$ corresponds to $F_k,...,F_{k+nN}$ such that the controller $u_i=F_i y_i$ steers any state at time $k$ to
zero at time $k+nN$.
\EPF%

In \REF{Proposition}{pro:specific-state-nullification} the system is assumed to be in a controller canonical form and the first coordinate of the
vectors $\{c_k\}_{k \in \Z}$ not zero. To prove \REF{Theorem}{thm:nullification} we need to show that $c_k \adj(A_k) b_k \ne 0$ if and only if the first
coordinate of $\tilde{c}_k$ (in the controller canonical form of the system) has a nonzero first coordinate. This fact is presented in
\REF{Proposition}{pro:c-tilde-vanish2} below. Towards this goal, we first prove that the property $c_k \adj(A_k) b_k \ne 0$ is an invariant of algebraic
equivalence.

\BCM \label{cla:c-tilde-vanish1} If \system and \tildesystem are algebraically equivalent by the transformation $\{T_k\}_{k \in \Z}$ then, for every $k
\in \Z$, $\det(T_{k+1})\tilde{c}_k\adj({\tilde{A}}_k)\tilde{b}_k=\det{(T_k)}c_k\adj(A_k)b_k$. \ECM
\BPF %
Define $\lambda_k = \frac{\det(T_{k+1})}{\det(T_k)}$ and $H^{}_k=\lambda^{}_k T^{}_k \adj(A^{ }_k)T^{-1}_{k+1}$. We
have,
\[\begin{aligned}
\tilde{A}_k H_k &= \lambda_k (T^{}_{k+1} A^{}_k T^{-1}_k) (T^{}_k \adj(A^{ }_k)T^{-1}_{k+1}) \\
                &= \lambda_k T^{}_{k+1} A^{}_k  \adj(A^{ }_k)T^{-1}_{k+1} \\
                &= \lambda_k T^{}_{k+1} \det(A^{ }_k)T^{-1}_{k+1} \\
                &= \lambda_k \det(A_k) I \\
                &= \det(T^{}_{k+1}A^{}_k T^{-1}_k) I \\
                &= \det(\tilde{A}_k) I.
\end{aligned}\]
Since the adjoint of a matrix is the only matrix such that $M \adj(M)=\det(M)I$, we get that
\[\adj(\tilde{A}_k)=\frac{\det(T_{k+1})}{\det(T_k)} T_k \adj(A_k)T_{k+1}^{-1}. \]
Therefore, \[ \tilde{c}_k \adj(\tilde{A}_k) \tilde{b}_k = \tilde{c}_k T_k^{-1} \adj(\tilde{A}_k) T_{k+1} \tilde{b}_k = \frac{\det(T_{k+1})}{\det(T_k)}
c_k \adj(A_k) b_k. \qedhere \]
\EPF %

\BP \label{pro:c-tilde-vanish2}If the system \tildesystem is in a controller canonical form~\eqref{equ:canonical-form}
and is algebraic equivalent to the system \system then, for every $k\!\in\!\Z$, the first coordinate of $\tilde{c}_k$
vanishes if and only if $c_k \adj(A_k) b_k =0$. \EP
\BPF %
The first entry of $\tilde{c}_k$ is given by $\tilde{c}_k (1,0,...,0)^T = \tilde{c}_k\adj(\tilde{A}_k)\tilde{b}_k$. By \REF{Claim}{cla:c-tilde-vanish1},
this is equal to $\frac{\det(T_{k+1})}{\det(T_k)}c_k\adj(A_k)b_k$.
\EPF %

Now we can conclude the proof of \REF{Theorem}{thm:nullification}.

\BPF[Proof of \REF{Theorem}{thm:nullification}] By \REF{Proposition}{pro:enough_to_nullify_a_given_state} and
\REF{Proposition}{pro:specific-state-nullification}, if the system is given in a controller canonical form and, for every $k \in \Z$, the first
coordinate of $c_k$ is not zero then it is memoryless output feedback nullifiable.

By \REF{Proposition}{pro:nullification-invariant}, it is enough to prove that the controller canonical form representation of the system is nullifiable.
By \REF{Proposition}{pro:c-tilde-vanish2}, the controller canonical form of a system satisfies the above condition if and only if  $c_k \adj(A_k) b_k
\ne 0$ for every $k \in \Z$. \EPF

\medskip

In the following theorem we provide an explicit bound on nullification time.

\BT \label{Thm:null-with-bound} If the system~\eqref{equ:system} is completely controllable, completely observable and $c_k \adj(A_k) b_k \ne 0$ for
every $k \in \Z$, then there exists a linear time-varying output-feedback controller of the form $u_k = F_k y_k$ that steers any initial state at any
time to the origin in $2(n^4+n^3+n^2)$ steps. \ET
\BPF%
By \REF{Claim}{cla:def_of_k_0}, we have that the index $k_0$ is smaller than $n^3+n^2$. For nullification we need an extra $n$ steps. Therefore, to
nullify a given initial state we need at most $n^3+n^2+n$ steps. To nullify any initial state we may need to repeat this procedure $n$ times (as
described in the proof of \REF{Proposition}{pro:enough_to_nullify_a_given_state}). Therefore, full nullification can be achieved in $n^4+n^3+n^2$ steps.
If the starting time is not fixed, we can apply the construction exposed in the proof of \REF{Proposition}{pro:specific-state-nullification} repeatedly,
at the cost of doubling the nullification time.
\EPF%

It is interesting to note that this bound, obtained for time-varying systems, is different than the bounds for time-invariant systems given
in~\citet{AW04}. This difference arise because when the vector $c$ is constant it is possible to use its properties to derive better bounds.

\section*{Acknowledgement}
I would like to thank my thesis advisor, Professor Zvi Artstein, for his support and patience. Without his guidance and insights, this research would
certainly not have been possible.

{
\bibliographystyle{ijc}
\bibliography{personalbib}
}

\end{document}